\title{ ~~\\ Non-Beiter ternary cyclotomic polynomials with an optimally large set of coefficients}
\author{Pieter Moree and Eugenia Rosu}
\documentclass[12pt]{article}
\usepackage{amssymb, latexsym, amsfonts, makeidx}
\usepackage{amsmath}
\def\@ptsize{2}
\newtheorem{Thm}{Theorem}
\newtheorem{Con}{Conjecture}
\newtheorem{Qu}{Question}
\newtheorem{Lem}{Lemma}
\newtheorem{Def}{Definition}
\newtheorem{cor}{Corollary}

\newtheorem{Prop}{Proposition}
\newcommand{\qed}{\hfill $\Box$}

\newcommand{\ds}{\displaystyle}
\begin{document}
\date{}
\maketitle 
{\def\thefootnote{}
\footnote{{\it Mathematics Subject Classification (2000)}. 
11N37, 11B83}}
\begin{abstract}
\noindent Let $l\ge 1$ be an arbitrary odd integer and $p,q$ and $r$ primes.
We show that there exist infinitely many ternary cyclotomic 
polynomials $\Phi_{pqr}(x)$ with $l^2+3l+5\le p<q<r$ such that the set of coefficients of each of them
consists of the $p$ integers in the interval $[-(p-l-2)/2,(p+l+2)/2]$. It is known that no larger coefficient range is possible. 
The Beiter conjecture states that the cyclotomic coefficients $a_{pqr}(k)$ of $\Phi_{pqr}$ satisfy
$|a_{pqr}(k)|\le (p+1)/2$ and thus the above family contradicts the Beiter conjecture. The two already known families of ternary
cyclotomic polynomials with an optimally large set of coefficients (found by G. Bachman) satisfy the Beiter conjecture.
\end{abstract}
\section{Introduction}
The $n$-th cyclotomic polynomial is defined by

$$\Phi_n(x)=\prod_{1\leq j\leq n, (j, n)=1}{(x-\zeta_n^j)}=\sum_{k=0}^{\varphi(n)}a_n(k)x^k,$$ with $\zeta_n=e^{2\pi i/n}$ and $\varphi(n)$ Euler's totient function. For $k>\varphi(n)$, take $a_n(k)=0$. We put $A\{n\}=\{a_n(k):k\ge 0\}$ and $A(n)=\max\{|a_n(k)|:k\ge 0\}$.
 
In this paper we will restrict to the (so-called ternary) case $n=pqr$, with $2<p<q<r$ primes. In \cite{Bachman-1} (Corollary 3) Bachman has proved that the difference between the largest and the smallest coefficient is at most $p$. Thus the cardinality 
of $A\{pqr\}$ is at most $p$. If it is exactly $p$ we say that $\Phi_{pqr}$ has an {\it optimally large set of coefficients} 

Sister Beiter \cite{Beiter} conjectured in 1968 that $|a_{pqr}(k)|\le (p+1)/2$. Beiter's conjecture implies that if $A\{pqr\}$ has
cardinality $p$, then
either $A\{pqr\}=\{-(p-1)/2,\ldots,(p+1)/2\}$ or $A\{pqr\}=\{-(p+1)/2,\ldots,(p-1)/2\}$. In 2004 Bachman \cite{Bachman} constructed two infinite
families of ternary cyclotomic polynomials having these respective coefficient sets.

Put $M(p)=\max\{A(pqr):r>q>p\}$, where $(q,r)$ runs over all prime pairs with $r>q>p$. Gallot and Moree \cite{GM1} disproved Beiter's conjecture. They showed that $M(p)>(p+1)/2$ for $p\ge 11$
and that, for $0<\epsilon<2/3$ one has $M(p)>(2/3-\epsilon)p$ for every $p$ sufficiently large. 
They conjectured that $M(p)\le 2p/3$. We like to point out that the construction of large coefficients given in this paper are
in accord with this conjecture.
Given an arbitrary integer $\delta\ge 0$ it
is thus conceivable that we can have {\it $\delta_{-}$-optimal}, respectively {\it $\delta_{+}$-optimal} ternary cyclotomic polynomials $\Phi_{pqr}$ having
$A\{pqr\}=\{-(p+1)/2+\delta,\ldots,(p-1)/2+\delta\}$, respectively $A\{pqr\}=\{-(p-1)/2-\delta,\ldots,(p+1)/2-\delta\}$ as coefficient sets.
Thus the construction of $0_{\pm}$-optimal $\Phi_{pqr}$ is due to Bachman. In this paper we give the construction (due to Rosu)
of $\delta_{\pm}$-optimal $\Phi_{pqr}$ for every $\delta\ge 1$ (with $\delta=(l+1)/2$). Every so constructed $\Phi_{pqr}$ will be
a counter-example to Sister Beiter's conjecture. 

The following lemma shows that it is enough to construct $\delta_{+}$-optimal $\Phi_{pqr}$.
\begin{Lem}
\label{aanvang}
Suppose that $\Phi_{pqr}$ is $\delta_{+}$-optimal and $r>pq$. Let $s>pq$ be
 a prime such that $s\equiv -r({\rm mod~}pq)$, then $\Phi_{pqs}$ is $\delta_{-}$-optimal.
\end{Lem}
The proof of this lemma is a corollary of the result of Kaplan \cite{Kaplan} that under the conditions of the
lemma one has $A\{pqr\}=-A\{pqs\}$. For a sharpening of this result, see Bachman \cite{Bachman1}.
\begin{Thm}
\label{main}
Let $l\geq 1$ be an odd integer and $p\ge l^2+3l+5$ a prime. Then there exists an infinite sequence of prime pairs 
$\{(q_j, r_j)\}_{j=1}^{\infty}$ with $pq_j<r_j$, $q_{j+1}>q_j$, such that 
$$\{a_{pq_jr_j}(k):k\geq 0\}=\Big\{-{p-l-2\over 2},\dots,{p+l+2\over 2}\Big\},$$ and hence $\Phi_{pq_jr_j}$ is 
$(l+1)/2_{+}$-optimal.
\end{Thm}
{\it Proof}. Theorem \ref{Main} (which is stated and proved in Section \ref{mainz}) and Lemma \ref{ineq}
allow us to determine $(q_j, r_j)$ such that
\begin{equation}
\label{falafel}
A\{pq_jr_j\}=\Big\{-{p-l-2\over 2},\dots,{p+l+2\over 2}\Big\}.
\end{equation}
(For notational convenience we write $-{p-1\over 2}$, rather than $-{(p-1)\over 2}$, etc..)
Using Dirichlet's theorem for primes in arithmetic progressions we see that there exists
an infinite family $\{(q_j,r_j)\}_{j=1}^{\infty}$ with $pq_j<r_j$, $q_{j+1}>q_j$, satisfying
(\ref{falafel}). \qed
\begin{cor} 
The Sister Beiter conjecture is false for every $p\ge 11$.
\end{cor}
Kaplan's Lemma (see Section \ref{kapell}) is the main tool in the proof of
Theorem \ref{Main}. Another helpful result we will use is the `jump-one property' (see \cite{BZ,GM2, Rosset}), 
which states that for ternary $n$ we have
$|a_n(k)-a_n(k+1)|\le 1$.
As a warm-up 
in using these tools the
reader may consult the alternative construction, by Gallot and Moree \cite{GM2}, of Bachman's $0_{\pm}$-optimal family. (When 
Bachman wrote his paper neither Lemma \ref{aanvang}, nor the jump-one property, nor 
the falsity of the Sister Beiter conjecture were known.)\\ 
\indent Let $M(p;q):=\max\{A(pqr):r>q\}$, where the maximum is over all primes $r$ exceeding $q$ and $p<q$ are fixed primes.
In Section \ref{harvesti} applications of the main result, Theorem \ref{Main} below, in the study of $M(p;q)$ will be given.
A question of Wilms (2010) first posed in \cite{GMW} will be answered in the positive. Also it is shown that the 
construction presented here yields a lower bound $M_R(p)$ for $M(p)$ that satisifies $M_R(p)\ge M_{GM}(p)$, with $M_{GM}(p)$ the lower
bound that was established in the paper by Gallot and Moree \cite{GM1}. Numerically one finds many $p$ with 
$M_R(p)>M_{GM}(p)$.\\ 
\indent The problem of finding non-Beiter ternary cyclotomic polynomials with an optimally large set of coefficients was
posed by Moree to Rosu as an internship problem. Rosu's solution of this is presented here. 
The write-up of this
construction was a cooperative effort of the authors, baving a rough draft by Rosu as a starting point. The section 
with the applications is due to the first author.
\section{Kaplan's Lemma}
\label{kapell}
Our main tool will be the following result of Kaplan \cite{Kaplan}, the
proof of which uses the identity
$$\Phi_{pqr}(x)=(1+x^{pq}+x^{2pq}+\cdots)(1+x+\cdots+x^{p-1}-x^q-\cdots-x^{q+p-1})
\Phi_{pq}(x^r).$$
\begin{Lem} {\rm (Nathan Kaplan, 2007)}.
\label{Kaplan}
Let $2<p<q<r$ be primes and $k\ge 0$ be an integer.
Put $$b_i=\begin{cases}a_{pq}(i) & \mbox{ if }ri\le k; \\ 0 & \mbox{ otherwise.}\end{cases}$$ We have
\begin{equation}
\label{lacheens}
a_{pqr}(k)=\sum_{m=0}^{p-1}(b_{f(m)}-b_{f(m+q)}),
\end{equation}
where $f(m)$ is the unique integer  such that
$f(m)\equiv r^{-1}(k-m) ({\rm mod~}pq)$ and $0\le  f(m) < pq$.
\end{Lem}
Lemma \ref{Kaplan} reduces the computation of $a_{pqr}(k)$ to that of $a_{pq}(i)$ for
various $i$. These binary cyclotomic polynomial coefficients are computed
in the following lemma. For a proof
see e.g. Lam and Leung \cite{LL}.
\begin{Lem}
\label{binary}
Let $p<q$ be odd primes. Let $\rho$ and $\sigma$ be the (unique) non-negative
integers for which 
\begin{equation}
\label{1pluspq}
1+pq=(\rho+1) p+(\sigma+1) q.
\end{equation}
Let $0\le m<pq$. Then either $m=\alpha_1p+\beta_1q$ or $m=\alpha_1p+\beta_1q-pq$
with $0\le \alpha_1\le q-1$ the unique integer such that $\alpha_1 p\equiv m({\rm mod~}q)$
and $0\le \beta_1\le p-1$ the unique integer such that $\beta_1 q\equiv m({\rm mod~}p)$.
The cyclotomic coefficient $a_{pq}(m)$ equals
$$\begin{cases}1 & \mbox{ if }m=\alpha_1p+\beta_1q \mbox{ with }0\le \alpha_1\le \rho,~0\le \beta_1\le
\sigma;\\ -1 & \mbox{ if }m=\alpha_1p+\beta_1q-pq \mbox{ with }\rho+1\le \alpha_1\le q-1,~\sigma+1\le 
\beta_1\le p-1;\\  0 & \mbox{ otherwise.}\end{cases}$$
\end{Lem}
We say that $[m]_p=\alpha_1$ is the {\it $p$-part of $m$} and $[m]_q=\beta_1$ is the {\it $q$-part
of $m$}. It is easy to see that
$$m=\begin{cases}[m]_pp+[m]_qq & \mbox{ if }[m]_p\le \rho \mbox{ and }[m]_q\le \sigma;\\ [m]_pp+[m]_qq-pq & \mbox{ if }[m]_p>\rho \mbox{ and }[m]_q>\sigma;\\ [m]_pp+[m]_qq-\delta_mpq & \mbox{ otherwise,}\end{cases}$$
with $\delta_m\in \{0,1\}$. Using this observation we find that, for $i<pq$,
$$b_i=\begin{cases}1 & \mbox{ if }[i]_p\le \rho, [i]_q\le \sigma \mbox{ and }[i]_pp+[i]_qq\le k/r;\\ -1 & \mbox{ if }[i]_p>\rho, [i]_q>\sigma \mbox{ and }[i]_pp+[i]_qq-pq\le k/r;\\ 0 & \mbox{ otherwise.}\end{cases}$$
Thus in order to evaluate $a_{pqr}(n)$ using Kaplan's lemma it suffices to compute $[f(m)]_p$, $[f(m)]_q$, 
$[f(m+q)]_p$ and $[f(m+q)]_q$. Indeed, as  $[f(m)]_p=[f(m+q)]_p$, it 
suffices to compute $[f(m)]_p$, $[f(m)]_q$, and $[f(m+q)]_q$. 
Note that, modulo $pq$, 
$$f(m+q)\equiv \frac{k-m}{r}-\frac{q}{r}\equiv f(m)+q\left(\left[-\frac{1}{r}\right]_p p+\left[-\frac{1}{r}\right]_q q\right)\equiv f(m)+q^2\left[-\frac{1}{r}\right]_q,$$
and hence
\begin{equation}
\label{q}
f(m+q)\equiv f(m)+q^2\left[-\frac{1}{r}\right]_q ({\rm mod~}pq).
\end{equation}
We will say that the $p$ and $q$-parts of $\gamma$ are in the same range if $0\leq[\gamma]_p\leq\rho$ and $0\leq[\gamma]_q\leq\sigma$ or if  $\rho+1\leq[\gamma]_p\leq q-1$ and $\sigma+1\leq[\gamma]_q\leq p-1$.

\section{The main construction}
\label{mainz}
Let $l\ge 1$ be an arbitrary odd integer.

The main construction only works if the residue $\ds 2/(l+2)$ modulo $p$ is in an union of two
intervals. The next lemma will be used to show that if $p\ge l^2+3l+5$, then the residue is in the
union of these two intervals.

\begin{Lem}
\label{ineq}
Let $(l+2)a\equiv 2 ({\rm mod~}p)$, where $0\leq a\leq p-1$ and $p\ge l^2+3l+5$. 
Then $a\in [l+2,(p-l-2)/2]\cup [(p+l+2)/2,p-l-2]$.
\end{Lem} 
\textit{Proof.}
The integer $a$ is of the form $(2+np)/(l+2)$ for some positive integer $n$ and hence $$\ds a\geq\frac{n(l+1)(l+2)+2}{l+2}> l+1.$$
Suppose that $$\ds p-1\geq \frac{np+2}{l+2}\geq p-l-1.$$ Then $$l+2-\frac{l+4}{p}\geq n \geq l+2-\frac{l^2+3l+4}{p}>l+1,$$ contradiction. Thus $a\leq p-l-2$.\\
\indent Finally suppose that $$\ds \frac{p-l}{2}\leq \frac{np+2}{l+2} \leq \frac{p+l}{2}.$$ Then 
$$\ds l+2-\frac{l(l+2)+4}{p}\leq 2n\leq l+2+\frac{l(l+2)-4}{p},$$ implying $2n=l+2$, which is not possible as $l$ is odd.\qed\\

The following lemma is important in the proof of the main construction:

\begin{Lem} \label{alpha} Let $u$ and $t$ be integers and $\alpha, \beta, r, p$ and $q$ be positive integers, with
$\beta>\alpha$ and $r>pq$.
Put $k(t)=ur+tpq$. 
Then 
\begin{equation}
\label{katee}
\alpha< \frac{k(t_1)}{r}<\beta,
\end{equation}
where $$t_1=\left[\frac{(\alpha-u)r}{pq}\right]+1.$$ 
\end{Lem}
\textit{Proof.} The inequality $\alpha < k(t)/r < \beta$ is equivalent 
with $$t\in \Big(\frac{(\alpha-u)r}{pq},\frac{(\beta-u)r}{pq}\Big).$$ 
Since this interval has length exceeding one, $t_t$ is in it and hence (\ref{katee}) is satisfied. \qed\\

It is practical to note that $1+pq=(\rho+1)p+(\sigma+1)q$ can be rewritten as $(p-1)(q-1)=\rho p+\sigma q$. By $[x]$ we
denote the largest integer $\le x$.
\begin{Thm}
\label{Main}
Let $l\ge 1$ be an arbitrary odd integer. Let $2<p<q$ be primes satisfying
$$q\geq \frac{(p+l)p}{2},~q\equiv \frac{2}{l+2} ({\rm mod~}p).$$ 
Let $\rho$ and $\sigma$ be the (unique) non-negative
integers for which $1+pq=(\rho+1) p+(\sigma+1) q$. 
Write $\rho=(p+l)s/2+\tau$ with $0\le \tau<(p+l)/2$ (thus $s=[2\rho/(p+l)]$).  
Let $r_1,r_2>pq$ be primes such that 
\begin{equation}
\label{r1r2}
-{1\over r_1}\equiv q-sp ({\rm mod~}pq){\rm ~and~}r_2\equiv -r_1 ({\rm mod~}pq).
\end{equation}
Let $\alpha_1^+=(p+l)q/2$ and $\alpha_1^-=(\rho +s)p-q$. 
Let $u_1^+,t_1^+,u_1^-,t_1^-,u_2^+,t_2^+,u_2^-,t_2^-$ be as in Table {\rm 1} and define 
$$k_1^+=r_1u_1^++t_1^+pq,~k_1^-=r_1u_1^-+t_1^-pq,~k_2^+=r_2u_2^+ +t_2^+pq,~k_2^-=r_2u_2^-+t_2^-pq.$$ 
Let $w(l+2)\equiv 2({\rm mod~}p)$, with $0\leq w\leq p-1$.\\ 
If $l+2\le w\leq (p-l-2)/2$, then
$$a_{pqr_1}(k_1^-)=-{p-l-2\over 2},~a_{pqr_1}(k_1^+)={p+l+2\over 2},$$
and $\{a_{pqr_1}(k):k\geq 0\}=\{-(p-l-2)/2,\dots,(p+l+2)/2\}$.\\
If  $(p+l+2)/2\le w\le p-l-2$, then 
$$a_{pqr_2}(k_2^-)=-{p-l-2\over 2},~a_{pqr_2}(k_2^+)={p+l+2\over 2},$$
and
$\{a_{pqr_2}(k):k\geq 0\}=\{-(p-l-2)/2,\dots,(p+l+2)/2\}$.\\
\indent In both cases we have $M(p;q):=\max\{A(pqr):r>q\}=(p+l+2)/2<2p/3$.
\end{Thm}
\begin{table}[ht]
\centerline{\bf TABLE 1}
\centering
\begin{tabular}{ |c | c | c |}
\hline
$u_1^+,t_1^+$ & $(\rho-\tau) p$ &  $\left[(\alpha_1^+-u_1^+)r_1/pq\right]+1$\\
$u_1^-,t_1^-$ & $-(p-l-2)q/2+((p-1)s+\tau)p$ &  $\left[(\alpha_1^--u_1^-)r_1/pq\right]+1$ \\
$u_2^+,t_2^+$ & $(p-1)q-(p-l-2)sp/2$ &  $\left[(\alpha_1^+-u_2^+)r_2/pq\right]+1$ \\
$u_2^-,t_2^-$ & $(p+l)q/2+\tau p$ &  $\left[(\alpha_1^--u_2^-)r_2/pq\right]+1$\\
\hline
\end{tabular}
\end{table}
The idea of the proof is the following:

\begin{itemize}
  \item In all cases, we take $k$ of the form $k=ut+pqr$, with $t=\left[(\alpha-u)r/pq\right]+1$. 
  \item We tabulate the $p$-parts and $q$-parts of $f(m)$, respectively of $f(m+q)$. We will claim what values each $b_{f(m)}$ must take, first for $0\leq m\leq p-1$, then for $q\leq m\leq q+p-1$.
  \item If the $p$-parts and $q$-parts of $f(m)$ are in different ranges, then $b_{f(m)}=0$ by Kaplan's Lemma and Lemma \ref{1pluspq}.
  \item If the $p$-parts and $q$-parts of $f(m)$ are in the same range we have to check whether $f(m)\le k/r$ or not, in
order to deduce that $b_{f(m)}=1$, respectively $b_{f(m)}=0$.

  \item In this way we find as a lower bound the number $\alpha$ that we chose in the definition of $k$ and as an upper bound a number $\beta$. 

  \item We check that $\alpha, \beta$ and $k$ verify the conditions of Lemma \ref{alpha}, which implies that 
    \begin{itemize}
      \item $k>0$, thus $a_{pqr}(k)$ is well-defined
      \item the inequalities claimed are true and the $b_{f(m)}$ take the claimed values.
      \item $a_{pqr}(k)$ can be computed by (\ref{lacheens})   
    \end{itemize}
\item We invoke the jump one property and a known upper bound on $M(p;q)$ to finish the proof.
\end{itemize}
We like to point out that in Kaplan's Lemma one is allowed to take any positive integer $k$. If applying it yields $a_{pqr}(k)\ne 0$,
then this implies that $0\le k\le \varphi(pqr)$. In particular $k_1^+,k_1^-,k_2^+,k_2^-$ are all in the range $(0,\varphi(pqr)]$.
An a posteriori check like this leaves one often with cleaner hands than checking this a priori.
\begin{Def}
We will use Kaplan's Lemma for the following $(k,r)$-pairs: $(k_1^{-},r_1)$, $(k_1^{+},r_1)$, $(k_2^{-},r_2)$, $(k_2^{+},r_2)$ and
denote the corresponding $f$-function by, respectively, $f_1^{-},f_1^{+},f_2^{-},f_2^{+}$.
\end{Def}
\begin{table}[ht]
\centerline{\bf TABLE 2}
\centering
\begin{tabular}{ |c c c c |}
\hline
$m$  & $[f_1^+(m)]_p$ & $[f_1^+(m)]_q$ & $[f_1^+(m+q)]_q$ \\
\hline
0 &  $\rho-\tau$ & 0 & $w$ \\
1 &  $\rho-\tau-s$  & 1 & $1+w$ \\
2 &   $\rho-\tau-2s$ & 2 & $2+w$ \\
 $\dots$ & $\dots$ & $\dots$ & $\dots$ \\
$\ds\frac{p+l}{2}-w$  & $\ds\rho-\tau-\left(\frac{p+l}{2}-w\right)s$ & $\ds\frac{p+l}{2}-w$ & $\ds\frac{p+l}{2}$ \\
%$\ds\frac{p+l}{2}-w+1$  & $\ds\rho-\tau-\left(\frac{p+l+2}{2}-w\right)s$ & $\ds\frac{p+l}{2}-w+1$ & $\ds\frac{p+l+2}{2}$ & 1 \\
$\dots$ & $\dots$ & $\dots$ & $\dots$ \\
$\ds\frac{p+l}{2}$ & 0 & $\ds\frac{p+l}{2}$ & $\ds\frac{p+l}{2}+w$ \\
$\ds\frac{p+l+2}{2}$ & $q-s$ & $\ds\frac{p+l+2}{2}$ & $\ds\frac{p+l+2}{2}+w$ \\
$\dots$ & $\dots$ & $\dots$ & $\dots$  \\
$p-1-w$  & $\ds q-\left(\frac{p-l-2}{2}-w\right)s$ & $p-1-w$ & $p-1$ \\
$p-w$  &$\ds q-\left(\frac{p-l}{2}-w\right)s$ & $p-w$ & 0 \\
$\dots$ & $\dots$ & $\dots$ & $\dots$ \\
$p-1$& $\ds q-\frac{p-l-2}{2}s$ & $p-1$ & $w-1$ \\
\hline
\end{tabular}
\end{table}
\noindent \textit{Proof of Theorem} \ref{Main}. In the case where $w=(p-l-2)/2$ or $w=(p+l+2)/2$ our argument
will need some minor modification and those cases will be left to the interested reader.\\
\indent If $l+2>(p-l-2)/2$ there is nothing to prove, so assume that $p\ge 3l+8$.
Note that $(p+l+2)q\equiv 2({\rm mod~}p)$. Since
$l\le p-4$, we infer that $q^*=(p+l+2)/2$, so $\sigma=q^*-1=(p+l)/2$. Knowing $\sigma$, we
compute $2\rho$ as 
\begin{equation}
\label{riedel}
\rho={pq-(l+2)q-2p+2\over 2p} 
\end{equation}
from (\ref{1pluspq}). We have $\rho<q/2$. Clearly we 
can write $\rho=(p+l)s/2+\tau$, with $0\leq \tau < (p+l)/2$.
Note that
$$q>{pq\over p+l}>{2\rho p\over p+l}\ge sp{\rm ~and~}q\ge {(p+l)p\over 2}\ge (\tau+1)p,$$
and hence
\begin{equation}
\label{maxstau}
q\ge \max\{sp+1,(\tau+1)p\}.
\end{equation}
\indent We have 
$$2\rho \ge {(p-(l+4))q\over p}\ge (p-(l+4)){(p+l)\over 2}\ge (l+2)(p+l),$$
and hence $s\ge l+2$.
Using that $1\le s< q$, we infer $(q-sp,pq)=1$. By
Dirichlet's theorem we then find that there are infinitely many primes
$r_1,r_2>pq$ satisfying congruence (\ref{r1r2}).\\

-{\tt The computation of $a_{pqr_1}(k_1^+)$}\\

\noindent {\tt Claim A}: Assume that $l+2\le w\le (p-l-4)/2$. Then Tables 2 and 3 are correct.\\
{\it Proof of Claim } A. Note that 
$f_1^+(0)\equiv \ds\frac{k_1^+}{r_1}\equiv u_1^+\equiv (\rho-\tau) p ({\rm mod~}pq)$.
By (\ref{r1r2}) we have $[-{1\over r_1}]_q=1$. Using (\ref{q}) and $q\equiv w({\rm mod~}p)$ we get 
\begin{equation}
\label{shift}
f_1^{+}(m+q)\equiv f_1^{+}(m)+q^2\left[-\frac{1}{r_1}\right]_q\equiv f_1^+(m)+wq ({\rm mod~}pq).
\end{equation}
This shows that the entries in the first row of Table 2 are as given. It is now very easy to compute
the remaining entries. We note that $f_1^+(m+1)\equiv f_1^+(m)-{1\over r_1}({\rm mod~}pq)$. Together with 
(\ref{r1r2}) it then follows that in order to get the $(j+1)$th row from the $j$th, we have to subtract $s$
from the $p$-part in row $j$  and add 1 to the $q$-part in row $j$ and reduce the result in such a way that 
the $p$-part is in $[0,q-1]$ and the $q$-part is in $[0,p-1]$. Since $\rho=(p+l)s/2+\tau$, we get at $p$-part
zero in row ${p+l\over 2}$. The correctness of Table 2 is established once we show that the last entry in
the second column is non-negative. In fact, we have
\begin{equation}
\label{rhorho}
q-\frac{p-l-2}{2}s> \rho.
\end{equation} 
Indeed, since $\rho=\frac{p+l}{2}s+\tau$ and $2\rho<q$, we have $\ds q-\rho>\rho>\frac{p-l-2}{2}s$. Using
(\ref{rhorho}) and $\sigma=\frac{p+l}{2}$ one eesily sees that the fourth and fifth column in Table 3
are correct.\\
\indent Next we consider $f_1^+(m)$ for the range $0\le m\le (p+l)/2$. Since $q>sp$ and
$f_1^+(m)$ for that range is increasing (we have
$f_1^+(m)=f_1^+(0)+m(q-sp))$, in order to show that  $f_1^+(m)\le k_1^+/r_1$ in
that range, we need only establish this inequality for 
$m=(p+l)/2$. Since $f_1^+(\frac{p+l}{2})=\frac{p+l}{2}q$ we have to
show that
\begin{equation}
\label{check1}
\frac{p+l}{2}q\leq\frac{k_1^+}{r_1},
\end{equation}
Recalling that $\alpha_1^+=(p+l)q/2$, $\ds k_1^+=u_1^+r_1+t_1^+pq$, 
with $\ds t_1^+=\left[(\alpha_1^+ -u_1^+)r_1/pq\right]+1$ and $pq<r_1$, 
it follows by Lemma \ref{alpha} that the inequality (\ref{check1}) 
indeed holds true.\\
\indent Since $\ds f_1^+(m) \geq f_1^+(\frac{p+l+2}{2})$ for 
$\ds\frac{p+l+2}{2}\leq m\leq p-1$, it is enough to check
\begin{equation}
\label{check2}
(q-s)p+\frac{p+l+2}{2}q-pq>\frac{k_1^+}{r_1},
\end{equation}
in order to verify the correctness of the sixth column in Table 3.
In order to establish (\ref{check2}), we will use Lemma \ref{alpha}. We see that
(\ref{check2}) follows
if $\beta-1\ge \alpha_1^+$, where $\beta$ denotes the left hand side in (\ref{check2}). Now 
since $\beta-\alpha_1^+=q-sp\ge 1$ this is obvious. The final column
in Table 3 is derived from the previous three, Lemma \ref{Kaplan} and Lemma \ref{binary}. \qed\\

\noindent {\tt Claim B}: Assume that $l+2\le w\le (p-l-4)/2$. Then Table 4 is correct.\\
{\it Proof of Claim} B. The second and third column are taken from Table 2, where
we used that $[f_1^+(m+q)]_p=[f_1^+(m)]_p$. To finish the proof we have to
show that
\begin{itemize}
 \item $[f_1^+(m)]_p\leq \rho$ and $[f_1^+(m+q)]_q\leq \sigma$ for $0\leq m \leq (p+l)/2-w$;
 \item $[f_1^+(m)]_p\leq \rho$ and $[f_1^+(m+q)]_q>\sigma$ for $(p+l+2)/2-w\leq m\leq (p+l)/2$; 
 \item $[f_1^+(m)]_p>\rho$ and $[f_1^+(m+q)]_q> \sigma$ for $(p+l+2)/2\leq m \leq p-1-w$;
 \item $[f_1^+(m)]_p>\rho$ and $[f_1^+(m+q)]_q\leq\sigma$ for $p-w\leq m\leq p-1$. 
\end{itemize}
In fact, these inequalities are obviously true, since $w\leq (p+l)/2=\sigma$ and (\ref{rhorho}) holds. \qed\\

\noindent {\tt Claim C}: We have $b_{f_1^+(m+q)}=0$ for $0\le m\le p-1$.\\
{\it Proof of Claim} C. {}From Table 4 we infer that 
$b_{f_1^+(m+q)}=0$ for $(p+l+2)/2-w\leq m\leq (p+l)/2$ and $p-w\leq m\leq p-1$ (since then the $p$-part and
$q$-part are not in the same range).
\begin{itemize}

  \item Since $f_1^+(m+q)\geq f_1^+(q)$ for $\ds 0\leq m\leq\frac{p+l}{2}-w$, it suffices to check that
\begin{equation}
\label{check3}
  (\rho-\tau) p+wq>\frac{k_1^+}{r_1},
\end{equation}  
in order to verify that $b_{f_1^+(m+q)}=0$ in this $m$-range.

  \item Since $\ds f_1^+(m+q)\geq f_1^+\left(q+\frac{p+l+2}{2}\right)$ for $\ds\frac{p+l+2}{2}\leq m\leq p-1-w$, it
suffices to check that
\begin{equation}
\label{check4}
  (q-s)p+\frac{p+l+2}{2}q+wq-pq>\frac{k_1^+}{r_1},
\end{equation}
 in order to verify that $b_{f_1^+(m+q)}=0$ in this $m$-range.
\end{itemize}
If both (\ref{check3}) and (\ref{check4}) are satisfied, then Claim C follows. By Lemma \ref{alpha} inequality
(\ref{check3}) follows if we can show that
$(\rho-\tau) p+wq-\frac{p+l}{2}q\ge 1$. Indeed, using (\ref{riedel}), (\ref{maxstau}) and $w\ge l+2$ (by assumption), we find
$$\ds(\rho-\tau) p+wq-\frac{p+l}{2}q=(w-(l+1))q-p+1-\tau p\geq 1+q-(\tau+1) p\geq 1.$$ 
That (\ref{check4}) is satisfied is obvious on noting that the left hand
side exceeds the left hand side in (\ref{check2}). \qed\\

By Kaplan's Lemma, Table 3 and Claim C we now infer that $a_{pqr_1}(k_1^+)=(p+l+2)/2$.

\begin{table}[ht]
\centerline{\bf TABLE 3}
\centering
\begin{tabular}{ |c | c | c | c | c | c | c|}
\hline
$m$                     & $[f_1^+(m)]_p$           & $[f_1^+(m)]_q$       & $[f_1^+(m)]_p$ & $[f_1^+(m)]_q$ & $f_1^+(m)$& $b_{f_1^+(m)}$\\
\hline
0                       &   $\rho-\tau$                 & 0                    & $\leq \rho$    &$\leq \sigma$ & $\leq k_1^+/r_1$ & 1\\
1                       &   $\rho-s-\tau$               & 1                    & $\leq \rho$    &$\leq \sigma$ &$\leq k_1^+/r_1$  & 1\\ 
2                       &   $\rho-2s-\tau$              & 2                    & $\leq \rho$    &$\leq \sigma$ & $\leq k_1^+/r_1$ & 1\\

$\dots$                 & $\dots$                  & $\dots$              & $\dots$& $\dots$& $\dots$& $\dots$\\

$\ds\frac{p+l}{2}$      & 0                        & $\ds\frac{p+l}{2}$   &$\leq \rho$     &$\leq \sigma$ &$\leq k_1^+/r_1$  & 1\\
$\ds\frac{p+l+2}{2}$    & $q-s$                    & $\ds\frac{p+l+2}{2}$ &$>\rho$    &$> \sigma$ &$> k_1^+/r_1$  & 0\\

$\dots$                 & $\dots$                  & $\dots$              & $\dots$& $\dots$& $\dots$& $\dots$\\

$p-1$                   &$\ds q-\frac{p-l-2}{2}s$  & $p-1$                &$>\rho$    &$> \sigma$ &$> k_1^+/r_1$  & 0\\

\hline
\end{tabular}
\end{table}
\noindent

\begin{table}[ht]
\centerline{\bf TABLE 4}
\centering
\begin{tabular}{| c |c| c |c |c |c| c | }
\hline
$m$                     & $[f_1^+(m+q)]_p$                            & $[f_1^+(m)]_q$        & $[f_1^+(m)]_p$ & $[f_1^+(m+q)]_q$\\
\hline
0                       & $\rho-\tau$                                    & $w$                     & $\leq \rho$  &$\leq\sigma$ \\

1                       & $\rho-\tau-s$                                  & $1+w$                   & $\leq \rho$  &$\leq\sigma$ \\

2                       & $\rho-\tau-2s$                                 & $2+w$                   & $\leq \rho$  &$\leq\sigma$\\

$\dots$                 & $\dots$                                   & $\dots$                 & $\dots$& $\dots$ \\

$\ds\frac{p+l}{2}-w$    & $\ds\rho-\tau-\left(\frac{p+l}{2}-w\right)s$   & $\ds\frac{p+l}{2}$      & $\leq \rho$  &$\leq\sigma$  \\

$\ds\frac{p+l+2}{2}-w$  & $\ds\rho-\tau-\left(\frac{p+l+2}{2}-w\right)s$ & $\ds\frac{p+l+2}{2}$    & $\leq \rho$  &$>\sigma$  \\

$\dots$                 & $\dots$                                   & $\dots$                 & $\dots$& $\dots$ \\

$\ds\frac{p+l}{2}$      & 0                                         & $\ds\frac{p+l}{2}+w$    & $\leq \rho$  &$>\sigma$ \\

$\ds\frac{p+l+2}{2}$    & $q-s$                                     & $\ds\frac{p+l+2}{2}+w$  & $>\rho$ &$>\sigma$  \\

$\dots$                 & $\dots$                                   & $\dots$                 & $\dots$& $\dots$ \\

$p-1-w$                 & $\ds q-\left(\frac{p-l-2}{2}-w\right)s$   & $p-1$                   & $>\rho$ &$>\sigma$  \\

$p-w$                   & $\ds q-\left(\frac{p-l}{2}-w\right)s$     & 0                       & $>\rho$ &$\leq\sigma$   \\

$\dots$                 & $\dots$                                   & $\dots$                 & $\dots$& $\dots$ \\

$p-1$                   &$\ds q-\frac{p-l-2}{2}s$                   & $w-1$                   & $>\rho$ &$\leq\sigma$    \\
\hline
\end{tabular}
Remark. For reasons of space we used the header $[f_1^+(m)]_p(=[f_1^+(m+q)]_p)$.
\end{table}
\noindent

-{\tt The computation of $a_{pqr_2}(k_2^+)$}\\

\noindent Here we assume that $(p+l+4)/2\le w\le p-l-2$, which ensures that
$l+2\le p-w\le (p-l-4)/2$.
Using Lemma \ref{alpha} one finds that $\alpha_1^+\leq \frac{k_2^+}{r_2}<\beta_2^+$, 
where 
$$\beta_2^+= \min\{(q-s)p+(p+l+2)q/2-pq, (\rho-\tau) p+(p-w)q\},~\alpha_1^+=(p+l)q/2.$$ By Lemma
\ref{alpha} it follows
from this that the analogues of (\ref{check1}), (\ref{check2}), (\ref{check3}) and (\ref{check4}) hold, where
we replace $k_1^+/r_1$ by $k_2^+/r_2$ and $w$ in (\ref{check3}) and (\ref{check4}) by $p-w$.\\ 
\indent We will show that
\begin{eqnarray}
\label{ktweeplus}
a_{pqr_2}(k_2^+) & = &\sum_{m=0}^{p-1}\left (b_{f_2^+(m)}-b_{f_2^+(m+q)}\right)
=\sum_{m=0}^{p-1}\left (b_{f_1^+(p-1-m)}-b_{f_1^+(p-1-m+q)}\right)\nonumber\\
&=& a_{pqr_1}(k_1^+)=\frac{p+l+2}{2}.
\end{eqnarray}
Note that $u_2^+\equiv (p-1)q-(p-l-2)sp/2\equiv u_1^+ +(p-1)(q-sp)({\rm mod~}pq).$ Using this observation it is easy to see that $\ds f_1^+(m)=f_2^+(p-1-m)$:
\begin{eqnarray*}
f_1^+(m) & \equiv & \frac{k_1^+}{r_1}-\frac{m}{r_1}\equiv u_1^+ +m(q-sp) \equiv u_2^++(p-1-m)(sp-q)\\
& \equiv & \frac{k_2^+}{r_2}-\frac{p-1-m}{r_2}\equiv f_2^+(p-1-m) ({\rm mod~}pq).
\end{eqnarray*} 
Since the analogues of (\ref{check1}) and (\ref{check2}) with $k_1^+/r_1$ replaced by $k_2^+/r_2$ hold, we
infer that $b_{f_1^+(m)}=b_{f_2^+(p-1-m)}$.\\
\indent Note that using (\ref{q}) we get $$f_2^+(m+q)\equiv f_2^+(m)+q^2\left[-\frac{1}{r_2}\right]_q\equiv f_1^+(p-1-m)+(p-w)q({\rm mod~}pq).$$  
By assumption we have 
$\ds l+2\leq p-w\leq (p-l-4)/2$.  As this is the only condition that we have used in obtaining the inequalities that involve $[f_1^+(m+q)]_q$, we can investigate whether the same bounds for $f_2^+(m+q)$ as for $f_1^+(p-1-m+q)$, 
where $0\leq m\leq p-1$, hold. Since as we have remarked (\ref{check3}) and (\ref{check4}) with 
$k_1^+/r_1$ replaced by $k_2^+/r_2$ and $w$ by $p-w$, hold, this is indeed so, and hence we conclude that 
$\ds b_{f_2^+(m+q)}=b_{f_1^+((p-1-m)+q)}$.
Thus we have
$$b_{f_2^+(m)}=b_{f_1^+(p-1-m)},~b_{f_2^+(m+q)}=b_{f_1^+(p-1-m+q)},~0\le m\le p-1,$$ 
implying (\ref{ktweeplus}).

\begin{table}[ht]
\centerline{\bf TABLE 5}
\centering
\begin{tabular}{|c c c c|}
\hline
$m$                      & $[f_1^-(m)]_p$                           & $[f_1^-(m)]_q$         & $[f_1^-(m+q)]_q$       \\
\hline
0                        & $\ds\rho+\frac{p-l-2}{2}s$               & $\ds\frac{p+l+2}{2}$   & $\ds\frac{p+l+2}{2}+w$ \\
1                        & $\ds\rho+\frac{p-l-4}{2}s$               & $\ds\frac{p+l+4}{2}$   & $\ds\frac{p+l+4}{2}+w$ \\
$\dots$                  & $\dots$                                  & $\dots$                & $\dots$                \\
$\ds\frac{p-l-4}{2}-w$   & $\ds\rho+(w+1)s$                         & $p-1-w$                & $p-1$                  \\
$\ds\frac{p-l-2}{2}-w$   &  $\ds\rho+ws$                            & $p-w$                  & 0                      \\
$\dots$                  & $\dots$                                  & $\dots$                & $\dots$                \\
$\ds\frac{p-l-4}{2}$     & $\rho+s$                                 & $p-1$                  & $w-1$                  \\
$\ds\frac{p-l-2}{2}$     & $\rho$                                   & $0$                    & $w$                    \\
$\dots$                  &$\dots$                                   & $\dots$                & $\dots$                \\
$p-w-1$                  & $\ds\rho-\left(\frac{p+l}{2}-w\right)s$  & $\ds\frac{p+l}{2}-w$   & $\ds\frac{p+l}{2}$     \\
$p-w$                    &$\ds\rho-\left(\frac{p+l+2}{2}-w\right)s$ & $\ds\frac{p+l+2}{2}-w$ & $\ds\frac{p+l+2}{2}$   \\
$\dots$                  & $\dots$                                  & $\dots$                & $\dots$                \\
$p-2$                    & $\tau+s$                                 & $\ds\frac{p+l-2}{2}$   & $\ds\frac{p+l-2}{2}+w$ \\
$p-1$                    & $\tau$                                   & $\ds\frac{p+l}{2}$     & $\ds\frac{p+l}{2}+w$   \\
\hline
\end{tabular}
Remark. If $w=\frac{p-l-4}{2}$ one should delete the first three rows in the table.
\end{table}

%$\ds k_1^-=u_1^-r_1+pqt_1^-$, where $\ds t_1^-=\left[\frac{(\alpha_1^- - u_1^-)r_1}{pq}\right]+1$ and $\ds u_1^-=(p-1)sp-\frac{p+l+2}{2}$, $\alpha_1^-=\rho p+sp-q$.

-{\tt The computation of $a_{pqr_1}(k_1^-)$}\\
 
\noindent {\tt Claim D}: Assume that $l+2\le w<(p-l-4)/2$. Then Tables 5 and 6 are correct. If $w=(p-l-4)/2$ they 
are also correct, but
with the first three rows (starting with $0,1,\ldots$) omitted in Table 5.\\ 
{\it Proof of Claim} D. 
Note that 
$$\frac{k_1^-}{r_1}\equiv \frac{p+l+2}{2}q+\left(\rho+\frac{p-l-2}{2}s\right)p-pq({\rm mod~}pq).$$  
Since $(p+l+2)/2<p$ and $\rho+(p-l-2)s/2<ps<q$ we have
$[f_1^-(0)]_p=(p+l+2)/2$ and $[f_1^-(0)]_q=\rho+(p-l-2)s/2$. This together with 
$f_1^{-}(m+q)\equiv f_1^-(m)+wq ({\rm mod~}pq)$, cf. the derivation of (\ref{shift}), yields the correctness of the
first row. {}From the first row, the remaining ones are easily determined, see the remarks made in the proof of
claim A.

\begin{table}[ht]
\centerline{\bf TABLE 6}
\centering
\begin{tabular}{|c | c | c | c | c | c | c|}
\hline
$m$                  & $[f_1^-(m)]_p$             & $[f_1^-(m)]_q$ & $[f_1^-(m)]_p$ & $[f_1^-(m)]_q$ & $f_1^-(m)$   &$b_{f_1^-(m)}$ \\
\hline
0                    & $\ds\rho+\frac{p-l-2}{2}s$ & $\ds\frac{p+l+2}{2}$  &$>\rho$ &$>\sigma$ &$\leq k_1^-/r_1$ & -1 \\
1                    & $\ds\rho+\frac{p-l-4}{2}s$ & $\ds\frac{p+l+4}{2}$  &$>\rho$ &$>\sigma$ &$\leq k_1^-/r_1$ & -1 \\

$\dots$              & $\dots$                    & $\dots$               &$>\rho$ &$>\sigma$ &$\leq k_1^-/r_1$ &-1  \\

$\ds\frac{p-l-4}{2}$ & $\rho+s$                   & $p-1$                 &$>\rho$ &$>\sigma$ &$\leq k_1^-/r_1$ & -1 \\
$\ds\frac{p-l-2}{2}$ & $\rho$                     & $0$                   &$\leq\rho$   &$\leq\sigma$   &$> k_1^-/r_1$ & 0 \\

$\dots$              & $\dots$                    & $\dots$               &$\leq\rho$   &$\leq\sigma$   &$>k_1^-/r_1$ &  0 \\

$p-2$                & $\tau+s$                   & $\ds\frac{p+l-2}{2}$  &$\leq\rho$   &$\leq\sigma$   &$> k_1^-/r_1$ & 0 \\
$p-1$                & $\tau$                     & $\ds\frac{p+l}{2}$    &$\leq\rho$   &$\leq\sigma$   &$> k_1^-/r_1$ & 0 \\
\hline
\end{tabular}
\end{table}
Note that
\begin{itemize}
  \item $[f_1^-(m)]_p>\rho$ and $[f_1^-(m)]_q>\sigma$ for $0\leq m \leq (p-l-4)/2$;
  \item $[f_1^-(m)]_p\leq\rho$ and $[f_1^-(m)]_q\leq\sigma$ for $(p-l-2)/2\leq m\leq p-1$,
\end{itemize}
showing the ocrrectness of columns four and five in Table 6.

\noindent Subclaim:\\
\noindent 1) $f_1^-(m)\le k_1^{-}/r_1$ for $\ds 0\leq m\leq (p-l-4)/2$;\\
2) $f_1^-(m)>k_1^{-}/r_1$ for $(p-l-2)/2\leq m\leq p-1$.\\
{\it Proof of Subclaim}. Recall that $sp<q$.
\begin{itemize}
  \item Since $\ds f_1^-(m) \leq f_1^- \left( \frac{p-l-4}{2} \right)$ for $\ds0\leq m\leq \frac{p-l-4}{2}$, we need to check that
  
  \begin{equation}
  \label{check5}
    (\rho+s)p+(p-1)q-pq=(\rho +s)p-q \leq \frac{k_1^-}{r_1},
  \end{equation}
in order to verify part 1.  
  \item Since $\ds f_1^-(m) \geq f_1^- \left( \frac{p-l-2}{2} \right)$ for $\ds\frac{p-l-2}{2}\leq m\leq p-1$, we have to check that
  
  \begin{equation}
  \label{check6}
    \rho p>\frac{k_1^-}{r_1},
  \end{equation}
in order to verify part 2.
\end{itemize}
Note that the inequality (\ref{check5}) is simply $\alpha_1^-\leq k_1^-/r_1$. 
We need to check that $\alpha_1^->0$. Since
\begin{equation}
\label{benhetzat}
\alpha_1^-\ge (\rho+1)p-q\ge (p-l-4)q/2>0, 
\end{equation}
this is easy.
In order to verify (\ref{check6}) we need
to show that $\rho p-\alpha_1^-\ge 1$, which is clear since 
$\rho p-\alpha_1^-=\rho p-(\rho+s)p+q=q-sp\ge 1$. This complets the proof of the sublclaim.\\ 
\indent Now the final column in Table 6 is deduced from the previous three
(by Kaplan's Lemma and Lemma \ref{binary}). \qed\\

\begin{table}[ht]
\centerline{\bf TABLE 7}
\centering
\begin{tabular}{| c | c | c | c | c | c |}
\hline
$m$                      & $[f_1^-(m+q)]_p$                           & $[f_1^-(m+q)]_q$ &  $[f_1^-(m)]_p$ & $[f_1^-(m+q)]_q$ \\
\hline
0                        & $\ds\rho+\frac{p-l-2}{2}s$               & $\ds\frac{p+l+2}{2}+w$ &$>\rho$ &$>\sigma$\\
1                        & $\ds\rho+\frac{p-l-4}{2}s$               & $\ds\frac{p+l+4}{2}+w$ &$>\rho$ &$>\sigma$\\
$\dots$                  & $\dots$                                  & $\dots$                & $\dots$ & $\dots$ \\
$\ds\frac{p-l-4}{2}-w$   & $\ds\rho+(w+1)s$                         & $p-1$                  &$>\rho$ &$>\sigma$\\
$\ds\frac{p-l-2}{2}-w$   & $\ds\rho+ws$                             & 0                      &$>\rho$ &$\leq\sigma$ \\
$\dots$                  & $\dots$                                  & $\dots$                & $\dots$& $\dots$ \\
$\ds\frac{p-l-4}{2}$     & $\rho+s$                                 & $w-1$                  &$>\rho$ &$\leq\sigma$ \\
$\ds\frac{p-l-2}{2}$     & $\rho$                                   & $w$                    &$\leq\rho$ &$\leq\sigma$ \\
$\dots$                  & $\dots$                                  & $\dots$                & $\dots$& $\dots$\\
$p-w-1$                  & $\ds\rho-\left(\frac{p+l}{2}-w\right)s$  & $\ds\frac{p+l}{2}$     &$\leq\rho$ &$\leq\sigma$ \\
$p-w$                    & $\ds\rho-\left(\frac{p+l+2}{2}-w\right)s$& $\ds\frac{p+l+2}{2}$   &$\leq\rho$ &$>\sigma$ \\
$\dots$                  & $\dots$                                  & $\dots$                & $\dots$& $\dots$\\
$p-2$                    & $\tau+s$                                 & $\ds\frac{p+l-2}{2}+w$ &$\leq\rho$ &$>\sigma$ \\
$p-1$                    & $\tau$                                   & $\ds\frac{p+l}{2}+w$   &$\leq\rho$ &$>\sigma$ \\
\hline
\end{tabular}
Remark. For reasons of space we used the header $[f_1^-(m)]_p(=[f_1^-(m+q)]_p)$.
\end{table}
\noindent Claim E: We have $b_{f_1^-(m+q)}=0$ for $0\leq m \leq p-1$.\\
{\it Proof of Claim} E. We assert that 
\begin{itemize}
 \item $[f_1^-(m)]_p> \rho$ and $[f_1^-(m+q)]_q> \sigma$ for $0\leq m \leq (p-l-4)/2-w$;
 \item $[f_1^-(m)]_p> \rho$ and $[f_1^-(m+q)]_q\leq\sigma$ for $(p-l-2)/2-w\leq m\leq (p-l-4)/2$;
 \item $[f_1^-(m)]_p\leq\rho$ and $[f_1^-(m+q)]_q\leq \sigma$ for $(p-l-2)/2\leq m \leq p-1-w$;
 \item $[f_1^-(m)]_p\leq\rho$ and $[f_1^-(m+q)]_q>\sigma$ for $p-w\leq m\leq p-1$. 
\end{itemize}
In fact, these are obviously true, since $w>0$ and $\tau\leq\rho$. We immediately infer that $b_{f_1^-(m+q)}=0$ for 
$(p-l-2)/2-w\leq m\leq (p-l-4)/2$ and $p-w\leq m\leq p-1$. 
  \begin{itemize}
   \item Since $f_1^-(m+q)\geq f_1^-(q)$ for $\ds 0\leq m \leq \frac{p-l-4}{2}-w$, we need to check that
      \begin{equation}
      \label{check7}  
      \left(\rho+\frac{p-l-2}{2}s \right)p+\left(\frac{p+l+2}{2}+w\right)q-pq>\frac{k_1^-}{r_1},
      \end{equation}
in order to show that $b_{f_1^-(m+q)}=0$ in this $m$-range.      
  \item Since $f_1^-(m+q)\geq f_1^-((p-l-2)/2+q) $ for $\ds \frac{p-l-2}{2}\leq m\leq p-w-1$, we need to check that
    \begin{equation}  
    \label{check8}      
    \rho p+wq>\frac{k_1^-}{r_1},
    \end{equation}
in order to show that $b_{f_1^-(m+q)}=0$ in this $m$-range.      
  
  \end{itemize}
If both (\ref{check7}) and (\ref{check8}) are satisfied, then 
Claim E follows.

We denote 
$$\beta_1^-=\min\{\rho p+\frac{p-l-2}{2}sp+\frac{p+l+2}{2}q+wq-pq, \rho p+wq \}.$$ Thus, in order to check 
inequalities (\ref{check7}) and (\ref{check8}) it is enough to show that $$\frac{k_1^-}{r_1}<\beta_1^-.$$
In order to prove this, we will use Lemma \ref{alpha}. Recalling that $\ds k_1^-=u_1^-r_1+t_1^-pq$, with 
$t_1^-=\left[(\alpha_1^- -u_1^-)r_1/pq\right]+1$ and $pq<r_1$, we need to check that 
$\beta_1^->\alpha_1^-$. For this we need to show that the following difference is strictly positive:
$$d:=\left(\rho+\frac{p-l-2}{2}s\right)p+\frac{p+l+2}{2}q+wq-pq-((\rho+s)p-q)$$
Using the assumption $w\geq l+2$, we get 
$d\ge d_1$ with $$2d_1=(p-l-4)sp+(3l+8-p)q.$$
Then, by first substituting $sp/2=(\rho p-\tau p)/(p+l)$, then using the inequality $q\ge (\tau+1)p$ and
finally using (\ref{benhetzat}), we obtain: 
\begin{eqnarray*}
d_1 & = & {p-l-4\over p+l}((\rho +1) p-(\tau +1) p) +{(3l+8)q-pq\over 2}\nonumber\\
& > & {p-l-4\over p+l}((\rho +1)p-q)+{(3l+8)q-pq\over 2}.\nonumber\\
& > & {p-l-4\over 2(p+l)}(p-l-4)q+{(3l+8)q-pq\over 2}.\nonumber
\end{eqnarray*}
Thus it suffices to show that 
$$d_2:={(p-l-4)^2\over p+l}+3l+8-p>0.$$
Simplifying the above expression gives
$$d_2={3l^2+8l+16\over p+l}>0.$$
Thus the conditions of Lemma \ref{alpha} are satisfied and hence
$0<\alpha_1^-\leq \frac{k_1^-}{r_1}<\beta_1^-$
and claim E is established. \qed\\

\noindent By Kaplan's Lemma, Table 6 and Claim E we now infer that 
\begin{equation}
\label{blomm}
a_{pqr_1}(k_1^-)=-\frac{p-l-2}{2}.
\end{equation}

-{\tt The computation of $a_{pqr_2}(k_2^-)$}\\
 
\noindent We claim that $\ds b_{f_2^-(m)}=b_{f_1^-(p-1-m)}$ and $\ds b_{f_2^-(m+q)}=b_{f_1^-(p-1-m+q)}$. This, using (\ref{blomm}) 
then yields $a_{pqr_2}(k_2^-)=a_{pqr_1}(k_1^-)=-(p-l-2)/2$, as required.

First, note that 
$$\ds \frac{k_2^-}{r_2}\equiv u_2^- \equiv \frac{p+l}{2}q+\tau p \equiv u_1^-+(p-1)(q-sp)({\rm mod~}pq).$$ 
Using this observation, it is easy to see that $f_2^{-}(m)=f_1^{-}(p-1-m_2)$:
\begin{eqnarray*}
f_2^-(m) & \equiv & \frac{k_2^- -m}{r_2}\equiv u_2^- +m(sp-q) \\
&\equiv & u_1^-+(p-1)(q-sp)+m(sp-q)\\
&\equiv & u_1^-+ (q-sp)(p-1-m)\\
& \equiv & \frac{k_1^- -(p-1-m)}{r_1}\equiv f_1^-(p-1-m)({\rm mod~}pq). 
\end{eqnarray*}

In order to have $\ds b_{f_2^-(m)}=b_{f_1^-(p-1-m)}$, we must have the same inequalities for $f_1^-(p-1-m)$ as for $\ds f_2^-(m)$.

Similarly, $\ds f_2^- (m+q)\equiv f_1^-(p-1-m)+(p-w)q$ (mod $pq$).  We want to have $\ds b_{f_2^-(m+q)}=b_{f_1^-(p-1-m+q)}$. 
The assumption we have made on $w$ ensures that $l+2\leq p-w\leq (p-l-4)/2$.  As this is the only condition that we have used in obtaining the inequalities that involve $[f_1^- (m+q)]_q$, we can ask for exactly the same bounds 
for $f_2^-(m+q)$ as for $f_1^- (p-1-m+q)$, where $0\leq m\leq p-1$.

Using Lemma \ref{alpha} one checks that the analogues of (\ref{check5}), (\ref{check6}), (\ref{check7}) and (\ref{check8}) with
$k_1^-/r_1$ replaced by $k_2^-/r_2$ and $w$ by $p-w$, hold true.
Hence we deduce that $b_{f_2^-(m)}=b_{f_1^-(p-1-m)}$ and $b_{f_2^-(m+q)}=b_{f_1^-(p-1-m+q)}$.\\

Thus we have finished the computation of $a_{pqr_1}(k_1^+)$, $a_{pqr_1}(k_2^+)$, $a_{pqr_2}(k_1^-)$ and
$a_{pqr_1}(k_2^-)$. Since the difference of the largest coefficient and the smallest coefficient in $\Phi_{pqr}(x)$
is at most $p$, cf.~\cite{Bachman}, we infer that 
$$\max A\{pqr_j\}=a_{pqr_j}(k_j^+),~\min A\{pqr_j\}=a_{pqr_j}(k_j^{-}).$$ 
By the jump one property, cf.~\cite{GM2}, we have
$|a_{pqr}(k)-a_{pqr}(k+1)|\le 1$ and hence we infer
that $A\{pqr_j\}=\{a_{pqr_j}(k_j^{-}),\ldots,a_{pqr_j}(k_j^+)\}$, as asserted. We have 
$M(p;q)\ge a_{pqr_j}(k_j^+)=(p+l+2)/2$. Recall that $q\equiv {2\over l+2}\equiv w({\rm mod~}p)$. We let $1\le w^*\le p-1$ be the 
inverse of $w$ modulo $p$. Note that $w^*=(p+l+2)/2$. We have $p/2<w^*< 3p/4$. By Lemma 3 of 
\cite{GMW} it then 
follows that $M(p;q)\le w^*=(p+l+2)/2$ and hence $M(p;q)=(p+l+2)/2$. 
Since $p\ge 3l+8$, it follows
that $(p+l+2)/2<2p/3$. \qed\\

\noindent {\tt Remark}. The upper bound on $M(p;q)$ can also be proved by directly invoking the main result
of B. Bzd{\c e}ga \cite{BZ}.
 
\section{Implications for the study of $M(p;q)$ and $M(p)$}
\label{harvesti}
Various papers in the literature study $M(p)$ and therefore it is perhaps natural to refine this to the study of
$M(p;q)$, as was first done by Gallot, Moree and Wilms \cite{GMW2}, see the MPIM-report \cite{GMW} for the full version. 
So far no algorithm is known to compute $M(p)$, whereas computing $M(p;q)$ is easy.

Gallot, Moree and Wilms \cite{GMW2} raise various questions, one of them (by Wilms), being:
\begin{Qu}
Given an integer $k\ge 1$, does there exist $p_0(k)$ and a function $q_k(p)$ such that if
$q\equiv 2/(2k+1)({\rm mod~}p)$, $q\ge q_k(p)$ and
$p\ge p_0(k)$, then $M(p;q)=(p+2k+1)/2$?
\end{Qu} 
Combining Theorem \ref{Main} with Lemma \ref{ineq} gives a positive
answer to this question (with $l=2k-1,~p_0(k)=4k^2+2k+3$ and $q_k(p)=(p+2k-1)p/2$).
\begin{Thm} 
\label{wilmspositive}
Let $l\geq 1$ be an odd integer and $p\ge l^2+3l+5$ a prime. Let $q\ge (p+l)p/2$ be a prime satisfying
$q\equiv 2/(l+2)({\rm mod~}p)$. Then
$M(p;q)=(p+l+2)/2$, where $M(p;q)=\max\{A(pqr):r>q\}$, where $r$ runs over all primes exceeding $q$.
\end{Thm}
Given a prime $p$ and a progression $a({\rm mod~}p)$ with $p\nmid a$, there are various results that 
ensure that
$M(p;q)$ assumes the same value, say $v$, for all primes $q\equiv a({\rm mod~}p)$ large enough. In this case we
write $m_p(a)=v$. 
Given an $1\le \beta\le p$ we let $\beta^*$ be the unique
integer $1\le \beta^*\le p-1$ with $\beta \beta^*\equiv 1({\rm mod~}p)$. \\
\indent Put
$${\cal B}_1(p)=\{\beta:~1\le \beta\le (p-3)/2,~\beta+\beta^*\ge p,~\beta^*\le 2\beta\},$$
$${\cal B}_{2}(p)=\{\beta:~1\le \beta\le (p-3)/2,~p\le \beta+2\beta^*+1,~\beta>\beta^*\},$$
$${\cal B}_{3}(p)=\{\beta:~1\le \beta\le (p-3)/2,~p\le 2\beta+\beta^*,~\beta\ge \beta^*\}.$$
$${\cal B}_{GM}(p):={\cal B}_{1}(p)\cup {\cal B}_{2}(p),~{\cal B}_R(p):={\cal B}_{1}(p)\cup {\cal B}_{3}(p).$$
Note that ${\cal B}_2(p)\subseteq {\cal B}_3(p)$ and hence ${\cal B}_{CM}(p)\subseteq {\cal B}_R(p)$. 
Furthermore, we have ${\cal B}_1(p)\cap {\cal B}_3(p)=\emptyset$.
Using 
Kloosterman sums and Weil's estimate, Cobeli \cite{Cobeli} can estimate the cardinality of these sets:
\begin{equation}
\label{bloeha}
\big|\#{\cal B}_{GM}(p)-\frac{p}{16}\big|\le 8\sqrt p(\log p +2)^3\,,
~~\big|\#{\cal B}_R(p)-\frac{p}{12}\big|\le 8\sqrt p(\log p +2)^3.
\end{equation}
The idea is that if say in the definition of ${\cal B}_1(p)$ we replace $\beta$ by $x$ and $\beta^*$ by $y$ and $x$ and $y$
are real numbers, we get a triangle $\Delta_1$. It is seen to have area asymptotic to $p^2/24$. Likewise, the triangles associated
to ${\cal B}_2(p)$ and ${\cal B}_3(p)$ have asymptotically area $p^2/48$, respectively $p^2/24$. Assuming now that the
inverses are uniformly distributed we expect asymptotically (as $p$ tends to
infinity) area$(\Delta)/p$ points $(\beta,\beta^*)$ inside the 
triangle $\Delta$. This indeed can be proved. Since $\Delta_1\cap \Delta_2=\emptyset$ and 
$\Delta_1\cap \Delta_3=\emptyset$, we then arrive at the main terms $p/16(=p/24+p/48)$ and 
$p/12=(p/24+p/24)$ in (\ref{bloeha}) above.\\
\indent The earlier large coefficient construction and the one presented in this paper can be formulated on
the same footing (as was pointed out to us by Yves Gallot).
\begin{Prop}$~$\\
\label{ziedubbel}
{\rm 1)} {\rm (Gallot and Moree)}. Let $\beta\in {\cal B}_{GM}(p)$. Then $m_p(\beta)\ge p-\beta $.\\
If $\beta\in {\cal B}_1(p)$ and $p=\beta+\beta^*$, then $m_p(\beta)=p-\beta $.\\
{\rm 2)} {\rm (Rosu)}. Let $\beta\in {\cal B}_R(p)$. Then $m_p(p-\beta^*)=p-\beta $.
\end{Prop}
{\it Proof}. 1) This is merely a consequence of the main theorem of Gallot and Moree \cite{GM1}.\\
2) We use the notation of Theorem \ref{Main},
Write $\beta=(p-l-2)/2$. Note that $\beta^*=p-w$. The inclusion $w\in [l+2,(p-l-2)/2]\cup [(p+l+2)/2,p-l-2]$
is seen to be equivalent with $\beta\in {\cal B}_R(p)$. Now invoke Theorem \ref{Main}.\qed\\

\noindent Let $S_i(p)$ be the set of $1\le v\le p$ for which part i of the latter proposition applies and yields
an identity or lower bound for $m_p(v)$ (thus $S_1(p)={\cal B}_{GM}(p)$ and 
$S_2(p)=\{p-\beta^*~:~\beta\in {\cal B}_R(p)\}$). A natural question that arises is to determine the intersection
$S_1(p)\cap S_2(p)$. Proposition \ref{vvorige} gives the answer, its proof depends on
the next proposition.
\begin{Prop}
Suppose that $\beta_1=p-\beta_2^*$ and $\beta_1,\beta_2\in {\cal B}_R(p)$.
Then $\beta_1,\beta_2\in {\cal B}_1(p)$ and $\beta_1=\beta_2$.
\end{Prop}
{\it Proof}. \\
-First case: $\beta_1,\beta_2\in {\cal B}_1(p)$.\\
The assumption $\beta_1=p-\beta_2^*$ implies
that $\beta_1^*=p-\beta_2$ and hence $\beta_1+\beta_1^*=2p-\beta_2-\beta_2^*$. Since
$\beta_i+\beta_i^*\ge p$, we infer that $\beta_i+\beta_i^*=p$. This together with $\beta_1=p-\beta_2^*$ yields
$\beta_1=\beta_2$.\\
-Second case: $\beta_1\in {\cal B}_1(p)$, $\beta_2\in {\cal B}_3(p)$.\\
{}From $\beta_1\in {\cal B}_1(p)$ and $\beta_1=p-\beta_2^*$ we infer that $\beta_1\ge \beta_1^*$. Now
$p\le \beta_1+\beta_1^*\le 2\beta_1$, contradicting $\beta_1\le (p-3)/2$.\\
-Third case: $\beta_1\in {\cal B}_3(p)$, $\beta_2\in {\cal B}_1(p)$.\\
We have $\beta_i\le (p-3)/2$ and hence $\beta_1+\beta_2\le p-3$. The assumption $\beta_1=p-\beta_2^*$ implies
that $\beta_1^*=p-\beta_2$. Since $\beta_1\in {\cal B}_3(p)$ we have $\beta_1^*\le \beta_1$ and hence
$p-\beta_2=\beta_1^* \le \beta_1$ and so $p\le \beta_1+\beta_2$. Contradiction.\\
-Fourth case: $\beta_1,\beta_2\in {\cal B}_3(p)$.\\
We have $\beta_i\le (p-3)/2$ and hence $\beta_1+\beta_2\le p-3$. Since $\beta_2\in {\cal B}_3(p)$, we have
$\beta_2^*\le \beta_2$ and thus from $\beta_1=p-\beta_2^*$ we infer that $p=\beta_1+\beta_2^*\le \beta_1+\beta_2$.
Contradiction. \qed
\begin{Prop}
\label{vvorige}
Let ${\cal P}_1$ be the set of primes $p\equiv 1({\rm mod~}4)$ such that the smallest solution, $x_0(p)$, of
$x^2+1\equiv 0({\rm mod~}p)$ satisfies $p/3\le x_0(p)\le (p-3)/2$. We have
$$S_1(p)\cap S_2(p)=
\begin{cases}
\emptyset & {\rm if~}p\not\in {\cal P}_1;\\
\{x_0(p)\} & {\rm otherwise}.
\end{cases}
$$
If $p\in {\cal P}_1$, then both parts of Proposition \ref{ziedubbel}
yield $m_p(x_0(p))=p-x_0(p)$.
\end{Prop}
{\it Proof}. By Proposition \ref{vvorige} we have 
$$S_1(p)\cap S_2(p)=\{\beta~:~1\le \beta \le (p-3)/2,~\beta+\beta^*=p,~\beta^*\le 2\beta\}.$$
The rest of the proof is left to the reader. \qed\\

\noindent This result together with Cobeli's estimate and Dirichlet's theorem for arithmetic progressions, shows that as 
$p$ tends to infinity, there is a set of primes $Q$ of density $\ge {7/48}$ such that $M(p;q)>(p+1)/2$. Thus
counter-examples to the Sister Beiter conjecture arise rather frequently.\\

\noindent Gallot et al.~\cite{GMW2} give conjectural values for $m_p(a)$ for $13\le p\le 23$.
Some of these can be shown to be true by Proposition \ref{ziedubbel}, which gives
$m_{13}(5)=8$, $m_{17}(12)=10$, $m_{19}(7)=11$, $m_{23}(16)=13$ and
$m_{23}(5)=14$. Of these only the first is not new, as $S_1(13)\cap S_2(13)=\{5\}$ and so also
part 1 of Proposition \ref{ziedubbel} yields $m_{13}(5)=8$.
\begin{Def}
Put 
$M_{GM}(p)=p-{\min}\{{\cal B}_{GM}(p)\}, M_{R}(p)=p-{\min}\{{\cal B}_R(p)\}$.
\end{Def}
Since ${\cal B}_2(p)\subseteq {\cal B}_3(p)$, we have
$$M(p)\ge M_R(p)\ge M_{GM}(p).$$
If $M_R(p)>M_{GM}(p)$, then the construction presented in this paper yields a better lower bound for $M(p)$ than
that established earlier by Gallot and Moree \cite{GM1}.
The primes $p<400$ with $M_R(p)>M_{GM}(p)$ are precisely:
$29,~37,~41,~83,~107,~109$, $149,~179,~181,~223,~227,~233,~241,~269,~281,~317,~347,~367,~379,~383,~389.$
\begin{Qu}
Are there infinitely many primes $p$ such that $M_R(p)>M_{GM}(p)$? If yes, give an estimate for the number
of such primes $\le x$.
\end{Qu}
Numerically it seems that with increasing $p$ occasionally larger and larger differences $M_R(p)-M_{GM}(p)$  occur.
\begin{Qu}
Is it true that $\lim \sup ~(M_R(p)-M_{GM}(p))=\infty$?
\end{Qu}

\section{Acknowledgement}
\indent Rosu would like to thank Pieter Moree at the Max Planck Institute for Mathematics in Bonn for suggesting the problem and for the great help and feedback he offered here during her two month-internship (summer 2010) at the MPIM Bonn.\\
\indent Moree likes to point out that Rosu was unaware of the question by Wilms' and that it hence came as a present surprise that
her construction as a by-product led to an answer of this question.\\
\indent The construction was found by paper and pencil work only and later underwent numerical checks by Yves Gallot. 
Furthermore, Gallot made many interesting remarks, in particular he introduced the set ${\cal B}_R(p)$. Merci
beaucoup, Yves!

\medskip\noindent {\footnotesize Max-Planck-Institut f\"ur Mathematik,\\
Vivatsgasse 7, D-53111 Bonn, Germany.\\
e-mail: {\tt moree@mpim-bonn.mpg.de}}\\

\medskip\noindent {\footnotesize Department of Mathematics, Evans Hall,\\
UC Berkeley, Berkeley CA 94720-3840, USA.\\
e-mail: {\tt rosu@math.berkeley.edu}}\\


{\small
\begin{thebibliography}{99}
\bibitem{Bachman-1} G. Bachman, On the coefficients of ternary cyclotomic 
polynomials, {\it J. Number Theory}  {\bf 100}  (2003), 104–-116.
\bibitem{Bachman} G. Bachman, Ternary cyclotomic polynomials with an optimally large set of coefficients, 
{\it Proc. Amer. Math. Soc.} {\bf 132} (2004), 1943--1950.
\bibitem{Bachman1} G. Bachman, On ternary inclusion-exclusion polynomials, {\it Integers  10}  (2010), A48, 623–-638. 
\bibitem{Beiter} M. Beiter, Magnitude of the coefficients of the cyclotomic 
polynomial $F\sb{pqr}\,(x)$, {\it Amer. Math. Monthly}  {\bf 75} (1968), 370--372.
\bibitem{BZ} B. Bzd{\c e}ga, Bounds on ternary cyclotomic coefficients, 
 {\it Acta Arith.}  {\bf 144}  (2010), 5--16. 
\bibitem{Cobeli} C. Cobeli, unpublished manuscript.
\bibitem{GM1} Y. Gallot  and P. Moree, Ternary cyclotomic polynomials having a large coefficient. 
{\it J. Reine Angew. Math.} {\bf 632} (2009), 105--125.
\bibitem{GM2} Y. Gallot  and P. Moree, Neighboring ternary cyclotomic coefficients differ by at most one, 
{\it J. Ramanujan Math. Soc.}  {\bf 24}  (2009),  235--248.
\bibitem{GMW} Y. Gallot, P. Moree and R. Wilms, The family of ternary cyclotomic polynomials with one free prime,
MPIM-preprint 2010-11, pp. 32.
\bibitem{GMW2} Y. Gallot, P. Moree and R. Wilms, The family of ternary cyclotomic polynomials with one free prime, {\it Involve}, to
appear, http://front.math.ucdavis.edu/1110.4590.
\bibitem{Kaplan} N. Kaplan, Flat cyclotomic polynomials of order three, 
{\it J. Number Theory} {\bf 127} (2007), 118--126. 
\bibitem{LL}T.Y. Lam and K.H. Leung, On the cyclotomic polynomial $\Phi_{pq}(X)$, {\it Amer.
Math. Monthly} {\bf 103} (1996), 562--564.
%\bibitem{Emma} E. Lehmer, On the magnitude of the coefficients of the cyclotomic polynomials, 
%{\it Bull. Amer. Math. Soc.} {\bf 42} (1936) 389--392.
%\bibitem{M} H. M\"oller, \"Uber die Koeffizienten des $n$-ten Kreisteilungspolynoms, 
%{\it Math. Z.} {\bf 119} (1971), 33--40. 
\bibitem{Rosset} S. Rosset, The coefficients of cyclotomic like polynomials of order 3, unpublished manuscript (2008) pp. 5. 
\end{thebibliography}}
\end{document}